\documentclass[12pt]{amsart}

\numberwithin{equation}{section}
\theoremstyle{plain}
\newtheorem{theorem}[equation]{Theorem}

\newtheorem{proposition}[equation]{Proposition}
\newtheorem{lemma}[equation]{Lemma}

\theoremstyle{remark}
\newtheorem{remark}[equation]{Remark}

\theoremstyle{definition}
\newtheorem{definition}[equation]{Definition}

\newcommand{\A}{{\mathcal A}}

\newcommand{\C}{{\mathcal C}}

\renewcommand{\L}{{\mathcal L}}

\newcommand{\Q}{\mathbb Q}
\newcommand{\R}{\mathbf{R}}

\newcommand{\V}{{\mathcal V}}

\newcommand{\lra}{\longrightarrow}
\newcommand{\ra}{\to}
\newcommand{\restr}{\mbox{\Large \(|\)\normalsize}}

\newcommand{\invlim}{\operatorname{\underset{\longleftarrow}{\lim}}}
\newcommand{\dirlim}{\operatorname{\underset{\longrightarrow}{\lim}}}

\newcommand{\Lip}{\operatorname{Lip}}

\newcommand{\loc}{\operatorname{loc}}

\renewcommand{\span}{\operatorname{span}}

\def\D{\partial}
\newcommand{\al}{\alpha}

\def\ga{\gamma}
\def\Ga{\Gamma}

\def\lra{\longrightarrow}

\def\ra{\to}

\def\si{\sigma}

\def\be{\beta}

\def\th{\theta}

\def\defeq{:=}

\newcommand{\no}{\noindent}

\newcommand{\mint}{-\!\!\!\!\!\!\int}
\newcommand{\ux}{\underline{x}}

\def\XXint#1#2#3{{\setbox0=\hbox{$#1{#2#3}{\int}$}
     \vcenter{\hbox{$#2#3$}}\kern-.5\wd0}}

\numberwithin{equation}{section}
\newcommand{\lla}{\longleftarrow}

\renewcommand{\hom}{\operatorname{Hom}}

\def\XXint#1#2#3{{\setbox0=\hbox{$#1{#2#3}{\int}$}
     \vcenter{\hbox{$#2#3$}}\kern-.5\wd0}}

\begin{document}
  \begin{abstract}
In this paper 
we prove the differentiability of 
Lipschitz maps
$X\to V$, where $X$ is a 
complete
metric measure space satisfying a doubling condition and a Poincar\'e
inequality, and $V$ denotes
a Banach space with the Radon Nikodym Property (RNP).
The proof depends on a new characterization of the differentiable
structure on such metric measure spaces, in terms of directional derivatives 
in the direction of 
tangent vectors to suitable rectifiable curves.
\end{abstract}

\title[Differentiability of RNP valued Lipschitz maps]
{Differentiability of Lipschitz maps from metric measure spaces
to Banach spaces with the Radon Nikodym property}

\date{\today}
\author{Jeff Cheeger}
\address{J.C.\,: Courant Institute of Mathematical Sciences\\
       251 Mercer Street\\
       New York, NY 10012}
\author{Bruce Kleiner}
\address{B.K.\,: Mathematics Department\\
         Yale University\\
             New Haven, CT 06520}
\thanks{The first author was partially supported by NSF Grant
DMS 0105128 and the second by NSF Grant DMS 0701515}
\maketitle
\tableofcontents

\section{Introduction}
\label{introduction}

In this paper we will use the term {\em PI space} to 
refer to a   $\lambda$-quasi-convex complete 
metric measure space  $(X,d^X,\mu)$
satisfying a doubling condition
\begin{equation}
\label{doubling}
\mu(B_{2r}(x))\leq 2^\kappa\cdot\mu(B_r(x))\,  ,
\end{equation}
and $p$-Poincar\'e inequality,
\begin{equation}
\label{poincare}
\mint_{B_r(x)}\,|f-f_{x,r}|\, d\mu
\leq \tau r\left(\mint_{B_{\lambda r}(x)}g^p\,d\mu\right)^{\frac{1}{p}}\, ,
\end{equation}
where 
$x\in X$, $r\in (0,\infty)$,  
$f$ is a continuous function,
$g$ is an upper gradient for $f$,
$$
\mint_A f\,d\mu\defeq \frac{1}{\mu(A)}\int_A f\,d\mu\, ,
$$
$$
f_{x,r}\defeq \mint_{B_r(x)}f\,d\mu\,; 
$$
see \cite{HeKo,cheeger,hei}.
We will also assume that the collection of measurable
sets is the completion of the Borel $\si$-algebra
with respect to $\mu$ (i.e. every subset of a set of 
measure zero is measurable) \cite[p.221]{royden}.
Sometimes, by abuse of 
language, we just say that $X$ is a PI space.  The notation above,
in particular the space $X$, the measure $\mu$, as well
as the constants $\kappa$ and $\lambda$,  
will be maintained 
throughout the paper.

In \cite{cheeger}, a differentiation theory for real valued Lipschitz
functions on PI spaces was given. The notion of differentiation is
expressed in terms of an atlas. An 
atlas consists of a countable collection
$\{(U_\alpha,y^\alpha)\}_{\al\in {\mathcal A}}$
of {\it charts},
where the $U_\alpha$'s are measurable subsets,
$\mu(X\setminus \bigcup_{\alpha\in {\mathcal A}} U_\alpha)=0$,
\mbox{} 
$y^\alpha:X\to \R^{k(\alpha)}$ is Lipschitz, and the charts satisfy
certain  additional conditions. 
We put $y^\alpha=(y^\alpha_1,\cdots,y_{k(\alpha)}^\alpha)$.

 Let $V$ denote
a Banach space.
\begin{definition}
\label{diffdef}
A  Lipschitz map
 $f:X\to V$ is differentiable
 almost everywhere  with respect to the
atlas $\{(U_\al,y^\al)\}_{\al\in {\mathcal A}}$  
 if
there is a  collection 
$$
\left\{\frac{\partial f}{\partial y^\alpha_m}:U_\al\lra V
\right\}_{\al\in {\mathcal A},\;1\leq m\leq k(\al)} 
$$ 
of Borel measurable functions
uniquely determined $\mu$-almost everywhere,
such that for almost every $\ux\in U_\al$,
\begin{equation}
\label{diff}
f(x)=f(\ux)+\sum_{m=1}^{k(\alpha)}
\frac{\partial f}{\partial y^\alpha_m}(\ux)
(y^\alpha_m(x)-y^\alpha_m(\ux)) + o(d^X(x,\ux))\, .
\end{equation}
\end{definition}
\no
We will  say that $f$ is differentiable at a specific
point $\ux\in X$ if (\ref{diff}) holds for that point.

The case $V=\R$ of Definition \ref{diffdef} was considered in \cite{cheeger};
one of the main results there \cite[Theorem 4.38]{cheeger} was
the existence of an atlas with respect to which 
every Lipschitz function
$f:X\to \R$ is differentiable almost everywhere. 
We will fix such an atlas  throughout the paper.
It follows readily from the definitions that if $(U_\al,y^\al)$
and $(\bar U_\al,\bar y^{\bar \al})$ are charts from two such
 atlases, then the matrix of partial derivatives 
$\frac{\D y^\al_m}{\D \bar y^{\bar \al}_{\bar m}}$
is defined and invertible almost everywhere in the overlap
$U_\al\cap \bar U_{\bar\al}$.
This also yields
 bi-Lipschitz invariant measurable tangent bundle
$TX$.

The main result of this paper is:
\begin{theorem}
\label{maindiff}
Every Lipschitz map from $X$ into a
Banach space with the Radon-Nikodym Property
is differentiable $\mu$-almost everywhere.
\end{theorem}
\vskip2mm

We recall that a Banach space $V$ has the {\it Radon-Nikodym Property}
(RNP) if every Lipschitz map $f:\R\to V$ is differentiable almost everywhere 
with respect to Lebesgue measure. Since $\R$ is an example of
a PI space, Theorem \ref{maindiff} is optimal in the sense that
the class of Banach space targets considered is maximal.

Just as in \cite{cheeger}, \cite{GFDA}, the differentiation 
theorem above imposes strong restrictions on PI spaces which
bi-Lipschitz embed in RNP targets, and may therefore
be used to deduce nonembedding theorems.
\begin{theorem}
\label{5.1}
If $X$ admits a bi-Lipschitz embedding in a Banach space with
the RNP, then for $\mu$-a.e. $x\in X$, every tangent cone at
$x$  is bi-Lipschitz homeomorphic to a Euclidean space.
\end{theorem}

\vskip2mm

\no
{\bf Discussion of the proof.}

The proof of Theorem \ref{maindiff} exploits the framework
introduced in \cite{GFDA}, which involves inverse systems of
finite dimensional Banach spaces and
their inverse limits; see Section \ref{il} for the relevant definitions.
It was observed in \cite{GFDA} that any separable Banach space $V$ can be realized
as a subspace of an inverse limit space $V\subset \invlim W_i$, where 
$\invlim W_i$  is the inverse limit of an inverse system of 
 finite dimensional Banach spaces. An advantage of this view point is that
the 
differentiability
 theory for real valued Lipschitz functions
 leads immediately to a natural notion of a {\it weak derivative} of 
Lipschitz map $f:X\to V$, which is a map taking values in $\invlim\, W_i$.

 It was pointed out in \cite{dprnp} that inverse limits $\invlim W_i$
of inverse systems of finite dimensional Banach spaces
are precisely the duals of separable Banach spaces.  Since $V$ has the RNP,
by a result of Ghoussoub-Maurey \cite{GhMau}, 
one can choose an embedding $V\subset \invlim W_i$ as above, so that 
the pair $(\invlim W_i,V)$ has 
the {\it Asymptotic Norming Property} (ANP). 
The ANP
was introduced by James-Ho
\cite{jamesho}, who showed that it implies the RNP,
  compare Section \ref{mm} and the appendix. 

The first step in the proof  of Theorem \ref{maindiff} is to 
show that if the weak derivative of 
$f$ takes values in the subspace $V\subset\invlim\, W_i$, then $f$ is
differentiable $\mu$-a.e.  The argument for this is brief, and
illustrates the
smooth interaction 
between the inverse limit setup,
 the ANP, and basic theorems of measure theory
(Egoroff's and Lusin's theorems).  As another illustration
of this smooth interaction, in the appendix we give a short
proof of the James-Ho theorem \cite{jamesho}.

The remainder of the proof, which appears in Section
\ref{secvelocitiesspan},
 is devoted to proving
that the weak derivative of $f$ takes values in $V$.  
A heuristic argument for this goes as follows.
If $c:I\to X$ is a Lipschitz curve, then  the composition
$f\circ c$ is differentiable almost everywhere because
$V$ has the RNP.  Hence the weak derivative of
$f\circ c$ coincides with its usual derivative, and in particular, 
lies in $V$.
When the curve $c$ has a well-defined measurable velocity 
$c':I\ra TX$
and the chain rule is applicable,  for almost 
every $t\in I$ the weak derivative evaluated
on $c'(t)$ will be the same as
$(f\circ c)'(t)$, which belongs to $ V$. 
 In this way one reduces the proof to showing that
for a full measure set of points $x\in X$, the tangent 
space $T_xX$ is 
spanned by the velocities of such curves $c$.  This fact, 
which is
of independent interest for the geometry of PI spaces, is established
in Section \ref{secvelocitiesspan}.

\bigskip

\no
{\bf Relation with previous work.}

By using an embedding $V\subset \invlim\, W_i$ as described above,
 a version of Theorem \ref{maindiff} was proved
  in \cite{GFDA} for a 
class of separable targets with a property that was
termed {\it Good Finite Dimensional Approximation} (GFDA). 
It was shown in \cite{GFDA} that separable dual spaces have the 
GFDA property.
An essential ingredient in
the proof of the differentiation theorem of \cite{GFDA} was
to show that if $V$ is a GFDA, then $V=\invlim \, W_i$.  It follows 
trivially
that for GFDA targets, the weak derivative lies in $V$; compare the 
discussion above.
 With the observation in \cite{dprnp} that inverse limits
$\invlim W_i$ are just the duals
of a separable Banach spaces, it followed that the class of GFDA's is precisely
the class of separable dual spaces, a strictly smaller class than that
of separable spaces with the RNP (or equivalently, the ANP);
see \cite{maccartneyobrien},  \cite{bourgaindelbaen}, \cite{bourgin}.
 Indeed, it was shown \cite{dprnp} that the GFDA condition is equivalent
to the ANP supplemented by an additional condition.

\section{Inverse limits and the ANP}
\label{il}

In this section 
we briefly recall some
facts from \cite{GFDA}, \cite{dprnp}.

Let $V$ denote a separable subspace of the dual space $Y^*$
of a separable Banach space $Y$. The pair $(Y^*,V)$ is said to have
the Asymptotic Norming Property (ANP) if for every sequence
$\{v_k\}\subset V$ the conditions
$$
v_k\stackrel{weak^*}{\longrightarrow} w\in Y^*\, ,
$$
$$
\|v_k\|\to \|w\|\, ,
$$ imply 
 {\it strong convergence},  i.e.
$\lim_{k\to\infty}\|v_k-w\|=0$. By \cite{GhMau1}, \cite{jamesho}, a separable
Banach space has the RNP if and only if it is isomorphic to some
$V$ 
which belongs to a pair 
$(Y^*,V)$ with the ANP.

We  observed  in \cite{dprnp} that a Banach space is
the dual of a separable Banach space if and only if
it is isometric to the inverse limit 
$\invlim\, W_i$ of an inverse inverse system,
\begin{equation}
\label{is}
W_1\stackrel{\th_1}{\lla} W_2\stackrel{\th_2}{\lla}\ldots
\stackrel{\th_{i-1}}{\lla} W_{i}\stackrel{\th_{i}}{\lla}\ldots\, ,
\end{equation}
where the $W_i$'s are finite dimensional Banach spaces
and the bonding maps $\theta_i$ are $1$-Lipschitz.
Such an inverse system will be called a {\it standard inverse system}.

 For the remainder of the paper we consider only standard inverse systems.

We recall that by definition, $\invlim\, W_i$ consists of all sequences
$(w_1,w_2,\ldots)$,
 where $w_i\in W_i$, 
 the 
compatibility condition  $\theta_i(w_{i+1})=w_i$ holds for
all $I$, and $\sup_i\,\|w_i\|<\infty$. The norm on $\invlim\, W_i$ is 
defined by 
$\|w\|=\lim_{i\to\infty} \|\pi_i(w)\|$. 

Let $\pi_j:\invlim\, W_i\to W_j$ denote the natural projection.
We may view the inverse limit $\invlim\,W_i$ as the dual space of the 
direct limit
$\dirlim\,W_i^*$ of the dual direct system
$$
W_1^*\stackrel{\th_1^*}{\ra}W_2^*\stackrel{\th_2^*}{\ra}\ldots\,
$$
  Then 
 a norm bounded sequence $\{v_k\}\subset\invlim\, W_i$ weak* converges
to $v_\infty\in \invlim_i\,W_i$ if and only if
 the  projected sequence
$\{\pi_j(v_k)\}\subset W_j$ converges to $\pi_j(v_\infty)$, for all $j$.


\section{Weak derivatives}
\label{mm}

Let $T^*X$ denote the
measurable cotangent bundle of the PI space $X$, and 
let $f:X\to V$ denote a Lipschitz map which is 
differentiable almost everywhere
in the sense
of Definition \ref{diffdef}.  The differential
$df$ is the  bounded measurable section of 
$T^*X\otimes V$
whose expression in the canonical trivialization of 
$T^*X\otimes V$
over $U_\alpha$ is 
\begin{equation}
\label{partials}
df=\left(
\frac{\partial f}{\partial y^\alpha_1},\ldots,
\frac{\partial f}{\partial y^\alpha_{k(\alpha)}}\right)\, ;
\end{equation}
compare (\ref{diff}). 

Recall that by definition, the tangent bundle $TX$ is the dual bundle of
the cotangent bundle $T^*X$. 
It will be convenient to work with the derivative $D_xf:TX\to V$,
which 
coincides with the differential $df$ under the identification
$\hom(TX,V)\simeq T^*X\otimes V$. 
Thus,
\begin{equation}
\label{dirdef}
D_xf(z)=\sum_{m=1}^{k(\alpha)}\frac{\partial f}{\partial y^\alpha_m}(x)z_m\, ,
\end{equation}
where 
$x\in U_\al$ and $z=\sum_m\,z_m\frac{\D}{\D y^\al_m}\in TX_x$.

Let  $V\subset\invlim\, W_i$ denote any closed 
linear subspace and assume $f:X\to V$ is Lipschitz.
Put $f_i=\pi_i\circ f$. 
For $\mu$-a.e. $x\in X$,
if  $v\in T_xX$,
the  collection of directional derivatives $\{D_xf_i(v)\}$ determines 
a norm bounded compatible sequence in the inverse system
$\{W_i\}$, and we thereby obtain 
 a {\it weak derivative} $\{D_xf_i\}:TX\to \invlim\, W_i$.
 The weak derivative is {\em weakly measurable}, in the sense
that its composition with $\pi_j:\invlim\,W_i\ra W_j$ is 
measurable, for all $j$.

\begin{remark}
\label{remweakisstrong} When a Lipschitz map 
$f:X\ra V\subset \invlim\,W_i$
is differentiable almost everywhere, 
the weak derivative is the true derivative, i.e.
for  $\mu$-a.e. $x\in X$, we have $D_x f=\{D_xf_i\}$ .   
This follows readily from the definitions.  In particular,
it follows that in this case the weak derivative is a 
Borel measurable mapping.
\end{remark}

Thus far in this section we have not explicitly invoked the 
Poincar\'e
inequality.  The next proposition, which is the converse of the
measurability statement in the preceding remark, will make 
use of it.

\bigskip
\begin{proposition}
\label{propweakstronglycontinuous}
Let $f:X\ra V$ be a Lipschitz map, where $V$ is an arbitrary
Banach space, and suppose the weak derivative 
$\{Df_i\}:TX\ra \invlim\,W_i$ 
is a Borel measurable mapping, with respect to the measurable
vector bundle structure on $TX$.  
 Then $f$ is differentiable almost everywhere.
\end{proposition}
\proof
Since $\{Df_i\}:TX\ra \invlim\,W_i$ is measurable, by Lusin's theorem, for almost
every $\ux\in X$ there is an $\al\in {\mathcal A}$ and a measurable
subset $A\subset U_\al$, such that $\ux\in A$ is a density point of 
$A$, and $\{\frac{\D f_i}{\D y^\al_m}\}:U_\al\ra \invlim\,W_i$ is continuous 
on $A$ for all $m\in\{1,\ldots,k(\al)\}$.

Let $\ell$ denote the function on the right-hand side of
(\ref{diff}), where the partial derivative $\frac{\D f}{\D y^\al_m}$
is replaced by 
$$
\{D_{\ux}f_i\}\left(\frac{\D}{\D y^\al_m}\right)
=\left\{\frac{\D f_i}{\D y^\al_m}(\ux)\right\}\,.
$$ 
Put $\ell_i=\pi_i\circ \ell$.
Then $D_x\ell$ is constant in the canonical local trivialization
of $TX$ on $U_\alpha$, and by the assumed
 continuity of $\{\frac{\D f_i}{\D y^\al_m}\}$ on $A$,  for all $x\in A$ and all $i$,
we have
$$
\lim_{x\to\ux}\|D_xf_i-D_x\ell_i\|=0\, ,
$$
where the convergence is uniform in  $i$.  Hence, by the 
Poincar\'e inequality
applied to the function $f_i-\ell_i$, and the fact that $f$ is 
Lipschitz,
the quantity
\begin{equation}
\label{strongdiff}
\sup_{x\in B_r(\ux)}\;\frac{\|f_i(x)-\ell_i(x)\|}{r}
\end{equation}
tends to zero as $r\ra 0$,   uniformly in $i$. Thus at $\ux$,
the weak derivative is a true derivative.

\qed

\vskip2mm
\noindent
{\bf The proof the main theorem, modulo showing that 
the weak derivative lies in $ V$.}

\noindent
For the remainder of the paper, we let
 $f:X\ra V$ be as in the statement of Theorem
\ref{maindiff}, and  $V\subset\invlim\,W_i$ be an embedding
such that the pair $(\invlim\,W_i,V)$ has the ANP, as in 
Section \ref{il}.

\begin{lemma}
\label{leminvsuffices}
Suppose the weak derivative
$\{D_xf_i\}$ takes values  in  $V\subset\invlim\, W_i$ 
for $\mu$-a.e. $x\in X$.  Then $f$ is differentiable almost
everywhere.
\end{lemma}
\proof  By Proposition \ref{propweakstronglycontinuous}, it suffices to show that 
the weak derivative is measurable.  To verify this, it suffices to check
that for every $\al$ and
every finite measure subset $A_0\subset U_\al$, there is a nearly
full measure subset $A\subset A_0$ where 
$\{\frac{\D f_i}{\D y^\al_m}\}:A\ra 
\invlim\,W_i$ is continuous.

Since  the maps,  $x\to \frac{\partial f_i}{\partial y^\alpha_m}$,
$m=1,\ldots,k(\alpha)$,  are measurable, by Lusin's theorem,
given any subset of finite measure, these functions are 
uniformly  continuous, for all $i$, off a subset of arbitrarily small measure.
By Egoroff's theorem and Lusin's theorem, the
same holds for $\|\{\frac{\partial f_i}{\partial y^\alpha_m}\}\|$.

Since 
$\{D_x f_i\}$ takes values in  $V$ for
$\mu$-a.e. $x$, for any finite measure subset of $X$,
we may invoke the Asymptotic Norming Property of
 $(\invlim\, W_i, V)$ 
 to conclude that 
$\{D_x f_i\}$ is strongly continuous off a set of 
arbitrarily small measure.

\qed

\section{Velocities of curves}
\label{secvelocitiesspan}

To show that the weak derivative $\{Df_i\}:TX\ra \invlim\,W_i$
takes values in 
$V$, and thereby  complete the proof of Theorem \ref{maindiff}, 
we will use 
directional derivatives along rectifiable curves, as indicated
in the introduction.  To formalize this,
we need to make precise the notion of the velocity of a  curve.

\vskip2mm
\noindent
{\bf Velocities and the chain rule.}

Let $X_0\subset X$ be a  full $\mu$-measure subset such that for 
every $x\in X_0$, if $x\in U_\al\cap U_\be$,
 then $x$ is a point of 
differentiability of $y^\al_m$ with respect to the chart $y^\be$,
for all $m\in \{1,\ldots,k(\al)\}$. 

\begin{definition}
\label{defvelocity}
If $c:I\ra X$ is a Lipschitz curve, $t\in I$ is a point of 
differentiability
of $y^\al\circ c$ for all $\al \in {\mathcal A}$, and $c(t)\in X_0$, 
then
the {\em velocity of $c$ at $t$} is defined to be the tangent vector 
$$
c'(t)=\sum_{m=1}^{k(\al)}\;(y^\al_m\circ c)'(t)
\frac{\D}{\D y^\al_m}\in T_{c(t)}X\,.
$$
\end{definition}
Note that this definition makes sense because of the choice 
of the set 
$X_0$.

With this definition, the chain rule becomes:
\begin{lemma}
\label{lemchainrule}
Suppose $c:I\ra X$ is a Lipschitz curve, and
the velocity vector $c'(t)\in TX$ is defined.  

1. For any Lipschitz function $u:X\ra \R$ which is differentiable 
with respect to the
atlas $\{(U_\al,y^\al)\}$ at $x$, 
the derivative $(u\circ c)'(t)$ is defined, and
$$
(u\circ c)'(t)=(D_{c(t)}u)(c'(t))\,.
$$

2. If 
$c(t)\in X$ is a point of weak differentiability of $f:X\ra V$
(i.e. $x$ is a point of differentiability of $f_i$ for all $i$),  
 and $t$ is a point of differentiability of $f\circ c$, 
then
we have the following chain rule relating the derivative
of $f\circ c$ and the weak derivative of $f$:
$$
(f\circ c)'(t)=\{(f_i\circ c)'(t)\}=\{(D_{c(t)}f_i)(c'(t))\}\,.
$$
In particular, the weak derivative of $f$ in the direction of the
velocity vector $c'(t)$ lies in $V\subset\invlim\,W_i$.
\end{lemma}
The proof is straightforward.

\vskip2mm
\noindent
{\bf Velocities span the tangent space.}

\noindent
Next we prove the following:

\begin{theorem}
\label{thmvelocitiesspan}
Fix a countable collection $\Phi=\{\phi_i:X\ra V_i\}$
of Lipschitz maps into RNP Banach spaces.
Let $\V_0$ be the collection of tangent vectors
$v\in TX$ such that there is a $1$-Lipschitz curve $c:I\ra X$
and $t\in I$, where 
$v=c'(t)$, and $t$ is a point of differentiability of $\phi_i\circ c$
for all $i$.
Then there is a  set $Z\subset X$ with $\mu(X\setminus Z)=0$, such
that 
 $\V_0\cap T_zX$ spans the fiber $T_zX$ at every point $z\in Z$.
\end{theorem}

\begin{remark}
\label{remdensedirections}
Elsewhere we will show that for a full measure set of
$x\in X$, there is a dense set of directions 
in  $\V_0\cap T_xX$, i.e.  the set
of rays in $T_xX$ which 
intersect $\V_0$ is dense in $T_xX$. However, we will not
need this finer result here.
\end{remark}

\noindent
{\em Proof of Theorem \ref{maindiff} using 
Lemma \ref{leminvsuffices} and Theorem \ref{thmvelocitiesspan}.}
Applying Theorem \ref{thmvelocitiesspan} with $\Phi=\{f:X\ra V\}$,
we obtain a full measure subset $Z\subset X$ as in the Theorem.
Let $W\subset Z$ be a full measure subset where $f$ is weakly
differentiable.  Then for every $x\in W$, and every 
$v\in \V_0\cap T_xX$, we have $\{(D_xf_i)(v)\}\in V$, by part
2 of  Lemma \ref{lemchainrule}.   Since $\V_0\cap T_xX$
spans $T_xX$, it follows that $\{(D_xf_i)(T_xX)\}\subset V$.
Hence $f$ is differentiable almost everywhere by 
Lemma \ref{leminvsuffices}.

\bigskip
\bigskip
We now turn to the proof of Theorem \ref{thmvelocitiesspan}.

Let $\V$ be the (fiberwise) span of $\V_0$, i.e. 
$\V\cap T_xX=\span(\V_0\cap T_xX)$.  

We begin with 
a preview of the argument.  We will first show that $\V$
defines a measurable sub-bundle of $TX$.
If there is a positive
measure set of points $x\in X$ where $\dim(\V\cap T_xX)<\dim T_xX$,
then by Lusin's theorem, we may pass to a positive 
measure subset $A$ of some $U_\al$ with the same property,
where in addition $\V$ lies
in a continuous codimension $1$
sub-bundle $E$ of the $k(\al)$-dimensional bundle $TX\restr_A$, 
where the continuity is defined 
with respect to the bundle chart induced by
$y^\al$.  If  $p$ is a density point of $A$, and $u$ 
is a linear combination of coordinates $y^\al_1,\ldots,y^\al_{k(\al)}$
whose derivative
at $p$ has kernel $E
\cap T_pX$, then one
finds that there
is an  upper gradient $\rho$ for $u$ such that
$$
\lim_{\stackrel{x\ra p,}{\,x\in A}}\;\rho(x)=0.
$$
This implies that the  average of $\rho$ over $B_r(p)$ tends to 
zero as $r\ra 0$, which contradicts the nondegeneracy of $u$.

\bigskip
We now give the details. Our first
step is:

\begin{lemma}
The sub-bundle $\V\subset TX$ is measurable.
\end{lemma}
Prior to proving the lemma, we 
 recall some facts about Suslin sets \cite{federer}.
A subset of a metric
space is {\em Suslin} if it is the continuous image
of a Borel subset of a complete separable metric space.
Suslin sets in a complete, $\si$-finite 
Borel regular  measure space such as $X$, are 
$\mu$-measurable; this is why we assumed that
the measure $\mu$ is complete.
Note that if $Z$ is a complete separable metric
space, then the image of a Suslin set $S\subset Z$ under
a Suslin measurable mapping $\tau:Z\ra X$ is also Suslin
(because the graph of $\tau$ is a Suslin subset of $Z\times X$). 

\proof
In brief, the proof is a straightforward application of
the  facts about Suslin sets recalled above.
  
It suffices to show that for each $\al$, the restriction
of $\V$ to $U_\al$ is measurable. We may assume that 
$U_\al$ is Borel measurable, and that it is contained
in the set $X_0$ defined before Definition \ref{defvelocity}. 

Let $\Ga$ denote the space of $1$-Lipschitz maps $c:[0,1]\ra X$,
equipped with the compact-open topology.  Then $\Ga$ is a 
complete, separable metric space, since $X$ is complete and
doubling.  Consider the collection $S_0\subset \Ga\times [0,1]$
of pairs $(c,t)$ such that $c(t)\in U_\al$, and the composition
$\phi_i\circ c:[0,1]\ra V_i$ is differentiable at $t$ for all $i$.
This is easily seen to be a Borel set.   Also, the map
$\si:S_0\ra \R^{k(\al)}$ which sends $(c,t)\in S_0$
to 
$$
(D_{c(t)}y^\al)(c'(t))=(y^\al\circ c)'(t)
\in \R^{k(\al)}
$$  
is Borel measurable.

For $j\in \{1,\ldots,k(\al)\}$, let
 $T_j$ be the set of points $x\in U_\al$
where the fiber $\V\cap T_xX$ has dimension $j$. 
We claim that
$T_j$ is a Suslin subset of $X$.  To see this, let
$S_1$ be the set of $j$-tuples
$((c_1,t_1),\ldots,(c_j,t_j))\in S_0^j$ such that
$c_1(t_1)=\ldots =c_j(t_j)$;
this is a closed subset of $S_0^j$.  Then let
$S_2\subset U_\al\times \wedge^j\,\R^{k(\al)}$
be the image of  $S_1$ under the
Borel map $S_1\ra U_\al\times \wedge^j\R^{k(\al)}$ which sends 
$((c_1,t_1),\ldots,(c_j,t_j))$ to 
$$
(c_1(t_1),\si((c_1,t_1))
\wedge\ldots\wedge\si((c_j,t_j))
\in U_\al\times \wedge^j\R^{k(\al)}\,.
$$
Then $\cup_{k\geq j}\;T_k$ -- the set of points where the
fiber of  $\V$ has dimension at least $j$ --
is the projection of $S_2\cap (U_\al\times 
(\wedge^j\R^{k(\al)}\setminus \{0\}))$ to $U_\al$, and is
therefore 
a Suslin set and $\mu$-measurable.   
It follows that $T_j$
is $\mu$-measurable for all $j\in \{1,\ldots,k(\al)\}$.

Let $\bar T_j\subset T_j$ be a full measure Borel subset of 
$T_j$.  Then $\cup_{j=1}^{k(\al)}\, \bar T_j$ has full measure 
in $U_\al$.  

Fix $j\in \{1,\ldots,k(\al)\}$.  Let $G(j,k(\al))$ denote the Grassman manifold of 
$j$-planes in $\R^{k(\al)}$.
Then there is a well-defined map $\ga_j:\bar T_j\ra G(j,k(\al))$
such that for every $x\in \bar T_j$, the fiber $\V\cap T_xX$ maps
under $D_xy^\al:T_xX\ra \R^{k(\al)}$ to the subspace $\ga_j(x)$
of $\R^{k(\al)}$.  To see that $\ga_j$ is measurable, pick
an open subset of $G(j,k(\al))$, and observe that
its inverse image in $\bar T_j$ is a Suslin set, using a
construction similar to the one above.

\qed

\bigskip
\bigskip
We now return to the proof of Theorem \ref{thmvelocitiesspan}.

Suppose $\dim(\V\cap T_xX)$ is strictly smaller than
$\dim T_x X$ for a positive measure set of points
$x\in X$.   Then for some index $\al\in {\mathcal A}$, there is
 a measurable 
 subset $A\subset U_\al$ with $\mu(A)\in (0,\infty)$,
where the strict inequality $\dim(\V\cap T_xX)<k(\al)$ holds. 
By Lusin's theorem, we may assume without loss of generality 
that $\V|A$ is contained in a codimension $1$ sub-bundle $E$ 
of $TX|A$,
where $E$ is a continuous sub-bundle relative to the bundle
charts given by $y^\al$, i.e. the fiber of $E$ at $x\in A$
is the kernel of $\psi(x)\circ D_xy^\al:T_xX\ra \R$, for some 
continuous
map
$\psi:A\ra (\R^{k(\al))^*}\setminus \{0\}$.  
It follows from the defining property of our  atlas 
$\{(U_\al,y^\al)\}$ -- specifically
the almost everywhere uniqueness of coefficients appearing in 
(\ref{diff}) --
 that we may  also assume
that for every $p\in A$, every nontrivial 
linear combination of the coordinate
functions $y^\al_m$  has nonzero pointwise upper Lipschitz constant
at $p$.

Choose $p\in A$, and put $\bar\psi=\psi(p)$.

\begin{lemma}
\label{lemzetabound}
There is a continuous function $\zeta:A\ra [0,\infty)$ such
that $\zeta(x)\ra 0$ as $x\ra p$, and 
$$
|D_x(\bar\psi\circ y^\al)(v)|\leq \zeta(x),
$$ 
for every  $v\in \V_0\cap T_xX$.
\end{lemma}
\proof
Let $c$ and $t$ be as in the definition of $\V_0$,
so that $c'(t)=v$.  Then
 $ D_x(\bar\psi\circ y^\al)(v)=\bar\psi((D_xy^\al)(v))$. 
Also, the vector $(D_xy^\al)(v)$ has uniformly bounded
norm since $y^\al$ is Lipschitz, and it lies 
in the hyperplane $\ker\psi(x)\subset \R^{k(\al)}$, which
approaches $\ker\bar\psi$ as $x\ra p$. The lemma follows.
\qed

\bigskip
Define a function $\rho:X\ra [0,\infty)$ by $\rho(x)=\zeta(x)$
if $x\in A$, and $\rho(x)=L$ otherwise, where $L$ is the Lipschitz
constant of $\bar\psi\circ y^\al$.  We claim  that $\rho$ is an 
upper gradient for $\bar\psi\circ y^\al$.   To see this, we need
only show that if $c:I\ra X$ is a $1$-Lipschitz curve, then
for almost every $t\in I$, we have 
$$
|\bar\psi\circ y^\al\circ c)'(t)|\leq \rho\circ c(t)\,.
$$
If $t\in I$ is such that $c(t)\not\in A$, then this obviously
holds, since $\bar\psi\circ y^\al\circ c$ is $L$-Lipschitz.  
If $t\in I$
is a point such that $c(t)\in A$ and
 the derivatives $(y^\al\circ c)'(t)$ and
$(\phi_i\circ c)'(t)$ are defined for all $i$, then $c'(t)$
is defined and lies in $\V_0$.  Therefore the chain rule applies,
and 
$$
|(\bar\psi\circ y^\al\circ c)'(t)|
=|\bar\psi((D_{c(t)}y^\al)(c'(t))|
\leq \zeta(c(t))=\rho\circ c(t)\,
$$
by Lemma \ref{lemzetabound}.
The remaining points $t\in I$ have measure zero.

Now let $p\in A$ be a density point of $A$.
Applying the Poincar\'e inequality to $\psi\circ y^\al$ on 
balls $B(p,r)$, using the fact that $p$ is a density point of 
$Z$ and $\psi\circ y^\al$ is Lipschitz, we conclude the pointwise
upper Lipschitz constant of $\psi\circ y^\al$ at $p$ is zero:
$$
\limsup_{r\ra 0}\;
\frac{\sup\;\left\{\psi\circ y^\al(x)-\psi\circ y^\al(p)\mid x\in B(p,r)
\right\}}
{r}=0.$$

This is a contradiction to the choice of $A$, which completes the
proof of the theorem.

\qed

\section{A new characterization of the
minimal upper gradient}

In this section will give a new characterization of the 
minimal upper gradient.  We then  apply this to give
a different proof of
 Theorem \ref{thmvelocitiesspan}.  It will also play a role
in the proof of the stronger version of Theorem
\ref{thmvelocitiesspan} alluded to in Remark 
\ref{remdensedirections}.

\bigskip
\noindent
{\bf Generalized upper gradients and minimal upper
gradients.}

Recall that if $u:X\ra \R$ is a Lipschitz function, then
a Borel measurable function $g:X\ra[0,\infty]$ is 
a {\em generalized upper gradient}  if there is a 
sequence of Lipschitz functions
 $u_k:X\ra \R$, and a  sequence  $g_k\in L^p_{\loc}(X)$
such that $g_k$ is an upper gradient for $u_k$ for all $k$,
and $u_k\stackrel{L^p_{\loc}}{\ra}u$,
$g_k\stackrel{L^p_{\loc}}{\ra}g$ \quad \cite[Section 2]{cheeger}.  
This is equivalent
to being a  {\em $p$-weak upper gradient}, 
i.e. satisfing
the usual upper gradient condition for all but a set of curves
of zero $p$-modulus
\cite{shanmugalingam,cheeger}.  It was shown in \cite{cheeger} that
every Lipschitz function $u:X\ra \R$ has a {\em minimal upper
gradient}, which is a generalized upper gradient $g_f:X\ra \R$
with the property that 
every generalized upper gradient $g$ satisfies
$g\geq g_f$ almost everywhere.

\bigskip
\noindent
{\bf Negligible sets.}

Let $\C$ denote the space of $1$-Lipchitz curves $c:[0,1]\to X$. 
With the
metric  $d(c_1,c_2)=\max_{t\in [0,1]}\; d^X(c_1(t),c_2(t))$, the space 
$\C$ is complete and separable, as is
 $\C\times [0,1]$ equipped  with the product metric.

\begin{definition}
A subset $N\subset \C\times [0,1]$ will be called 
{\it negligible} if it is Borel, and
for all $c\in \C$,
$$
\L(N\cap (\{c\}\times [0,1]))=0\, ,
$$
where $\L$ denotes Lebesgue measure on $[0,1]$.
Clearly, a countable union of negligible sets is negligible.
\end{definition}

\begin{theorem}
\label{thmgcharacterization}
Let $f:X\ra \R$ be a Lipschitz function with minimal upper
gradient $g_f$, and let
$N\subset \C\times[0,1]$ be a negligible set. 
Define a function $\hat g_f:X\ra [0,\infty)$ by letting 
$\hat g_f(x)$ be
the supremum of the set 
$$
\{\;|(f\circ c)'(t)|\;\mid\;
(c,t)\in \C\times[0,1]\,\setminus\, N,\;\;c(t)=x,\;
\mbox{and}\; (f\circ c)'(t)\;\mbox{exists}\}\,
$$
if it is nonempty, and $0$ otherwise.
  Then $\hat g_f$ coincides with $g_f$ almost everywhere.
\end{theorem}

\proof
We begin
by showing that the function $\hat g_f$ defined above is
$\mu$-measurable.  For this it suffices to show
that $\hat g_f^{-1}((a,\infty))$ is a Suslin set for all 
$a\in [0,\infty)$.

Let $\pi:\C\times [0,1]\ra X$ be the map
$\pi((c,t))=c(t)$.   Then $\pi$ is continuous,
and $\hat g_f^{-1}((a,\infty))$ is the image under
$\pi$ of the Borel set
$$
\left\{(c,t)\in \C\times[0,1]\setminus N\;\mid\; 
(f\circ c)'(t)
\;\mbox{exists},\; 
|(f\circ c)'(t)|\in (a,\infty)\right\}\,;
$$
hence $\hat g_f^{-1}((a,\infty))$ is Suslin, as claimed.

Let $g:X\ra \R$ be a Borel measurable function such
that $g\geq \hat g_f$, and $ g=\hat g_f$ almost everywhere.
Notice that $ g$ is an upper gradient for $f$, because
if $c\in \C$, then for a.e. $t\in [0,1]$, the derivative
$(f\circ c)'(t)$ exists,  $(c,t)\not\in N$, and
$$
|(f\circ c)'(t)|\leq  \hat g_f\circ c(t)\leq g\circ c(t)\,.
$$
Therefore $\hat g_f\geq g_f$ almost everywhere, since
$g_f$ is a minimal upper gradient.

Observe that $\hat g_f\leq \Lip f$  everywhere; this follows
from the fact that if $c\in \C$ and $(f\circ c)'(t)$ exists,
then
$\Lip_{c(t)}(f)\geq |(f\circ c)'(t)|$ because $c$ is $1$-Lipschitz.
By \cite[Thm 6.1]{cheeger}, we have $g_f=\Lip f$ almost everywhere, 
and hence  $\hat g_f\leq g_f$ almost everywhere.  Thus
$\hat g_f=g_f$ almost everywhere.
\qed

\begin{remark}
The full strength of Theorem \ref{thmgcharacterization} is not 
used in the 
application given below.  It would be sufficient to
know that  $g\geq C\,g_f$
$\mu$-a.e., for some $C\in (0,\infty)$ which does not depend
on $f$, and this is considerably easier to prove, 
see \cite[Prop. 4.26]{cheeger}.
\end{remark}

\bigskip
\noindent
{\bf An alternate proof of Theorem \ref{thmvelocitiesspan}.}

Returning to the setting of Theorem \ref{thmvelocitiesspan},
let $N\subset\C\times[0,1]$
be the negligible set of pairs $(c,t)$ such that one
of the compositions
$\{y^\al\circ c\}_{\al\in{\A}}$ or 
$\{\phi_i\circ c\}_{\phi_i\in\Phi}$ is not differentiable at $t$.

For each $\al\in \A$, and each 
rational $k(\al)$-tuple $(a_1,\ldots,a_{k(\al)})\in \Q^{k(\al)}$,
let $h_{(a_1,\ldots,a_{k(\al)})}:X\ra \R$ be the function
$\hat g_f$
defined as
   in Theorem \ref{thmgcharacterization}, with 
$f=\sum_m\;a_m\,y^\al_m$.

Now let $Z\subset X$ be a full measure Borel set
such that:
\begin{enumerate}
\item $Z\subset X_0$, where $X_0\subset X$ is the 
set defined before Definition \ref{defvelocity}.
\item For all $\al\in\A$, $(a_1,\ldots,a_{k(\al)})\in \Q^{k(\al)}$,
and $x\in Z$,  
the function $h_{(a_1,\ldots,a_{k(\al)})}(x)$
 coincides with
the  pointwise upper Lipschitz constant of the function
$f=\sum_m\;a_m\,y^\al_m$ at every point in $Z$.
\item The function $h_{(a_1,\ldots,a_{k(\al)})}$ is approximately
continuous at every point in $Z$.
\end{enumerate}

Suppose that the $\V_0\cap T_xX$ does not span the fiber $T_xX$
at some point $x\in U_\al\cap Z$.  Then there is a nonzero
$k(\al)$-tuple $(b_1,\ldots,b_{k(\al)})\in \R^{k(\al)}$ such
that $\sum_m\,b_m\,(D_xy^\al_m)(v)=0$ for all $v\in \V_0\cap T_xX$.
The chain rule  (Lemma \ref{lemchainrule}) implies that $\hat g_u(x)=0$, where
$u\defeq \sum_m\,b_m\,y^\al_m$.  If $\{a^j\}\subset \Q^{k(\al)}$
is a sequence converging to $b=(b_1,\ldots,b_{k(\al)})$, then $h_{a^j}(x)\ra 0$.
Therefore  the upper pointwise Lipschitz constant of 
$\sum_m\,a_m^j\,y^\al_m$ at $x$ tends to zero as well.
Hence the nontrivial linear combination $\sum_m\,b_m\,y^\al_m$
has zero pointwise upper Lipschitz constant at $x$, which
contradicts the fact that the differentials of the $y^\al_m$'s
are independent on $U_\al\cap Z$.  

\qed

\section*{Appendix.  The ANP implies the RNP
(An alternate  proof of the theorem of James-Ho \cite{jamesho})}

Using the formalism of inverse limit spaces, in this
appendix we will give a 
short direct 
proof that the ANP implies
that a Lipschitz function $f:I\to V\subset\invlim W_i$,
 is differentiable a.e., provided the separable space $V$  has the ANP;
compare \cite{jamesho}.

Let $I\subset \R$ denote a finite interval.
Given $\eta>0$, by Lusin's theorem (respectively, Egoroff's and Lusin's 
theorems)
there exists
$A\subset I$ such that $\L(I\setminus A)<\eta$ and in addition, 
$D_xf_i$ (for all $i$) 
and 
$\|\{D_xf_i\}\|$ are uniformly continuous on $A$. 
It suffices to show that
$f$ is differentiable, with derivative $D_{\ux}f=\{D_{\ux} f_i\}$, for every 
density point
$\ux$ of $A$.

Put 
$$
\Delta_\epsilon f(\ux)=\frac{f(\ux+\epsilon)-f(\ux)}{\epsilon}\, .
$$
Then for all $i$,
$$
\pi_i(\Delta_\epsilon f(\ux))=\frac{f_i(\ux+\epsilon)-f_i(\ux)}{\epsilon}\, ,
$$
\begin{equation}
\label{1}
\pi_i(\Delta_\epsilon f(\ux))=\frac{1}{\epsilon}\cdot
\int_{\ux}^{\ux+\epsilon}D_xf_i \, dx\, .
\end{equation}
Taking the norm of both sides of (\ref{1}) and using 
$\|D_xf_i\|\leq \|\{D_xf_i\}\|$, gives
\begin{equation}
\label{2}
\|\pi_i(\Delta_\epsilon f(\ux))\|\leq \frac{1}{\epsilon}\cdot
\int_{\ux}^{\ux+\epsilon}\|\{Df_i \}\|\, dx\, .
\end{equation}

Letting $\epsilon\to 0$, then $i\to\infty$  in (\ref{1}) and using
that $D_xf_i$ is continuous on $A$ and that $\ux$ is a density point of $A$, 
shows that
$\Delta_\epsilon f$ converges weak${}^*$ to $\{D_{\ux} f_i\}$ 
as $\epsilon\to 0$.
Letting $i\to\infty$, then $\epsilon\to 0$ in (\ref{2}) and using
that $\|\{D_xf_i\}\|$ is continuous on $A$ and $\ux$ is a density point of $A$, shows that 
$\limsup_{i\to\infty}\|\Delta_\epsilon f\|\leq \|\{D_{\ux} f_i\}\|$.  Since
$\Delta_\epsilon f(\ux)\in V$, the ANP implies that $\Delta_\epsilon f(\ux)$
converges strongly to $\{D_{\ux} f_i\}$.

\bibliography{pirnp}
\bibliographystyle{alpha}
\addcontentsline{toc}{subsection}{References}

\end{document}